\documentclass[12pt,a4paper,reqno]{amsart}
\usepackage{amsmath}
\usepackage{amsfonts}
\usepackage{amssymb}
\usepackage{bm} 

\newcommand\eps{\varepsilon}
\newcommand\R{{\mathbf{R}}}
\newcommand\C{{\mathbf{C}}}

\renewcommand\L{{\mathcal{L}}}

\newcommand\loc{{\operatorname{loc}}}

\renewcommand\mod{{\ \operatorname{mod}\ }}
\newcommand\pc{{\operatorname{pc}}}

\theoremstyle{plain}
  \newtheorem{theorem}[subsection]{Theorem}

  \newtheorem{corollary}[subsection]{Corollary}

\theoremstyle{remark}
  \newtheorem{remark}[subsection]{Remark}

\theoremstyle{definition}

\include{psfig}

\begin{document}

\title[Pseudoconformal compactification of NLS]{A pseudoconformal compactification of the nonlinear Schr\"odinger equation and applications}
\author{Terence Tao}
\address{Department of Mathematics, UCLA, Los Angeles CA 90095-1555}
\email{tao@@math.ucla.edu}

\vspace{-0.3in}
\begin{abstract}
We interpret the lens transformation (a variant of the pseudoconformal transformation) as a pseudoconformal compactification of spacetime, which converts the nonlinear Schr\"odinger equation (NLS) without potential with a nonlinear Schr\"odinger equation with attractive harmonic potential.  We then discuss how several existing results about NLS can be placed in this compactified setting, thus offering a new perspective to view this theory.
\end{abstract}

\maketitle

\section{Introduction}

Let $d \geq 1$ be an integer.  We consider solutions $u: I \times \R^d \to \C$ of the
free (i.e. zero-potential) non-linear Schr\"odinger equation
\begin{equation}\label{nls}
(i \partial_t + \frac{1}{2} \Delta) u = \mu |u|^{p-1} u
\end{equation}
on a (possibly infinite) time interval $I$, where $p > 1$ is an exponent and $\mu = \pm 1$.  For the algebraic manipulations
below we shall assume that our solution has sufficient regularity and decay to justify all the formal calculations; standard limiting arguments (see e.g. \cite{caz}) then allow us to utilise the same computations for any regularity class for which one has a strong local wellposedness theory.

The case $\mu=+1$ is \emph{defocusing}, while the case $\mu=-1$ is \emph{focusing}.  This equation 
enjoys the scale-invariance
\begin{equation}\label{scaling}
u(t,x) \mapsto \frac{1}{\lambda^{2/(p-1)}} u(\frac{t}{\lambda^2}, \frac{x}{\lambda})
\end{equation}
for $\lambda > 0$, and also has the conserved \emph{mass}
$$ M(u) = M(u(t)) := \int_{\R^d} |u(t,x)|^2\ dx.$$
The \emph{$L^2$-critical} or \emph{pseudoconformal} power $p = 1 + \frac{4}{d}$ is of special interest for a number of reasons.  Firstly, the conserved mass $M(u)$ becomes invariant under the scaling \eqref{scaling}.  Secondly, one also acquires an additional invariance, the \emph{pseudoconformal invariance} $u \mapsto u_\pc$, where $u_\pc: -I^{-1} \times \R^d \to \C$
is defined on the time interval $-I^{-1} := \{-1/t: t \in I\}$ (assuming $0 \not \in I$) by the formula
\begin{equation}\label{pcdef}
u_\pc(t,x) := \frac{1}{|t|^{d/2}} u( \frac{-1}{t}, \frac{x}{t} ) e^{i|x|^2 / 2t}.
\end{equation}
The pseudoconformal transform is an involution (thus $(u_\pc)_\pc = u$) and preserves the conserved mass
$$ M(u_\pc) = M(u)$$
and (as observed in \cite{blue}) more generally preserves the $L^q_t L^r_x$ Strichartz norms for all admissible $(q,r)$ (i.e. $2 \leq q, r \leq \infty$ and $\frac{2}{q} + \frac{d}{r}=\frac{d}{2}$ with $(d,q,r) \neq (2,2,\infty)$):
$$ \| u_\pc \|_{L^q_t L^r_x(-I^{-1} \times \R^d)} =  \| u \|_{L^q_t L^r_x(I \times \R^d)}.$$
In particular we have
$$ \| u_\pc \|_{L^{2(d+2)/d}_{t,x}(-I^{-1} \times \R^d)} =  \| u \|_{L^{2(d+2)/d}_{t,x}(I \times \R^d)}.$$
We remark that if $p$ is not the pseudoconformal power, and $u$ solves \eqref{nls}, 
then $u_\pc$ does not solve \eqref{nls}, but instead solves the very similar equation
$$ (i \partial_t + \frac{1}{2} \Delta) u_\pc = \mu t^{\frac{d}{2}(p-1) - 2} |u_\pc|^{p-1} u_\pc.$$
Note that the pseudoconformal transformation inverts the time variable, sending $t=0$ to $t=\pm \infty$ and vice versa.
Because of this fact, this transform has been very useful in understanding the asymptotic behaviour of free nonlinear Schr\"odinger equations. 

The purpose of this note is to highlight a close cousin of the pseudoconformal transformation, namely the \emph{lens transform}
which we will define shortly.  This transform, introduced in the study of NLS in \cite{lens}, \cite{lens2}, \cite{carles0},
\emph{compactifies} the time interval to $(-\pi/2,\pi/2)$ rather than inverting it, and is thus very analogous to the
\emph{conformal compactification} map of Penrose \cite{penrose}, which has been useful in studying the asymptotic behaviour of nonlinear wave equations (see e.g. \cite{christ}).   The one ``catch'' is that the lens transform introduces an attractive quadratic potential $\frac{1}{2} |x|^2$ to the linear component of the Schr\"odinger equation, which changes the long-time dynamics (for instance, the propagators for the linear equation are now time-periodic with period $2\pi$).  The lens
transform was used in \cite{carles}, \cite{carles2} to study nonlinear Schr\"odinger equations with harmonic potential, but
we argue here that it can also be used to clarify much of the theory concerning the \emph{free} nonlinear Schr\"odinger equation, particularly the portion of the theory concerning scattering and uniform spacetime bounds.  (Similar ideas appeared in \cite{cw-rap}.)  In particular,
asymptotics at $t \to \pm \infty$ are converted to asymptotics at $t \to \pm \pi/2$, thus converting the global-in-time theory to local-in-time theory.

The results we present here are not new, being essentially due to earlier work by other authors; thus the paper here is more of a survey than a research paper.  However we believe that the unifying perspective afforded by the lens transform is not widely known, and thus hopefully of interest to readers.

The author thanks Tonci Cmaric for pointing out that the results of Begout and Vargas imply an inverse Strichartz theorem, and Remi Carles and Jim Colliander for helpful corrections, references, and comments.  The author also thanks the anonymous referee for helpful comments and corrections.

\section{The lens transform}

Given any function $u: I \times \R^d \to \C$, we define the \emph{lens transform} $\L u: \tan^{-1}(I) \times \R^d \to \C$ of $u$ on the time interval
$\tan^{-1}(I) := \{ \tan^{-1}(t): t \in I \} \subset (-\pi/2,\pi/2)$, where $\tan^{-1}: \R \to (-\pi/2,\pi/2)$ is the 
arctangent function, by the formula
$$ \L u(t,x) := \frac{1}{\cos^{d/2} t} u( \tan t, \frac{x}{\cos t} ) e^{-i |x|^2 \tan t / 2}.$$
Thus for instance, $\L u(0,x) = u(0,x)$, or in other words the lens transform does not distort the initial data.  Its inverse
is given by
$$ \L^{-1} v(t,x) = \frac{1}{(1+t^2)^{d/4}} v( \tan^{-1} t, \frac{x}{\sqrt{1+t^2}} ) e^{i |x|^2 t / 2(1+t^2)}.$$
Like the pseudoconformal transform, the lens transform also preserves the mass and Strichartz norms:
$$ M(\L u) = M(u); \quad \| \L u \|_{L^q_t L^r_x( \tan^{-1}(I) \times \R^d )} = \| u \|_{L^q_t L^r_x(I \times \R^d)}.$$
In particular the map $u(\tan t) \mapsto \L u(t)$ is unitary in $L^2_x(\R^d)$ for each $t$, and the pure Strichartz norm is preserved: 
$$ \| \L u \|_{L^{2(d+2)/d}_x( \tan^{-1}(I) \times \R^d )} = \| u \|_{L^{2(d+2)/d}_x(I \times \R^d)}.$$
As observed in \cite{lens}, \cite{lens2}, \cite{carles0}, if $u$ solves \eqref{nls}, then $\L u$  
solves a nonlinear Schr\"odinger equation with attractive harmonic potential:
\begin{equation}\label{lens}
 (i \partial_t + \frac{1}{2} \Delta - \frac{1}{2} |x|^2) \L u = \mu |\cos t|^{\frac{d}{2}(p-1) - 2} |\L u|^{p-1} \L u.
\end{equation}
This is especially useful in the pseudoconformal case $p = 1 + \frac{4}{d}$, in which case we simply have
\begin{equation}\label{lu-form}
(i \partial_t + \frac{1}{2} \Delta - \frac{1}{2} |x|^2) \L u = \mu |\L u|^{4/d} \L u.
\end{equation}
More generally, the equation \eqref{lens} can be sensibly extended to all times $t \in \R$ provided that we are in the scattering-subcritical regime $p > 1 + \frac{2}{d}$, so that the weight $|\cos t|^{\frac{d}{2}(p-1) - 2}$ is locally integrable.  In the case of the scattering-critical power $p = 1 + \frac{2}{d}$ or the scattering-supercritical powers $p < 1 + \frac{2}{d}$ the equation \eqref{lens} has more serious singularities at $t = \pm \pi/2$ and the asymptotics are more nonlinear here; see \cite{ozawa}.

The lens transform is closely related to the pseudoconformal transformation, indeed one easily verifies the formula
$$ \L u_\pc(t,x) = \L u( t + \pi/2 \mod \pi, x ) \hbox{ for } t \in (-\pi/2,\pi/2)$$
where $s \mod \pi$ is the unique translate of $s$ by an integer multiple of $\pi$ which lies in the fundamental domain $(-\pi/2,\pi/2]$; note that $\tan( t + \pi/2 \mod \pi ) = -1/\tan(t)$.  
Thus the lens transform conjugates the pseudoconformal transformation to (essentially) 
a translation in the lens time variable (which is the arctangent of the original time variable).  Because of this, many arguments in the literature which rely on the pseudoconformal transformation can easily be recast using the lens transform instead.  To give one very simple example, observe that the \emph{harmonic energy} of $\L u(t)$,
\begin{equation}\label{harmonic}
 \int_{\R^d} \frac{1}{2} |\nabla \L u(t,x)|^2 + \frac{1}{2} |x|^2 |\L u(t,x)|^2 + \frac{\mu}{p+1} |\cos t|^{\frac{d}{2}(p-1) - 2} |\L u(t,x)|^{p+1}\ dx
 \end{equation}
is equal to the \emph{classical energy} of $u(\tan t)$,
$$ \int_{\R^d} \frac{1}{2} |\nabla u(\tan t,x)|^2 + \frac{\mu}{p+1} |u( \tan t, x)|^{p+1}\ dx$$
plus the \emph{pseudoconformal energy} of $u(\tan t)$,
$$ \int_{\R^d} \frac{1}{2} |(x + i\tan t\nabla) u(\tan t,x)|^2 + \frac{\mu \tan^2 t}{p+1} |u( \tan t, x)|^{p+1}\ dx.$$
In the pseudoconformal case $p = 1 + \frac{4}{d}$, the harmonic energy of $\L u$ is conserved.  Since the classical energy of $u$ is also conserved, we conclude the conservation law for the pseudoconformal energy.  For other values of $p$, the harmonic energy of $\L u$ enjoys a monotonicity formula, which yields the standard monotonicity formula for the pseudoconformal energy.
The scattering space $\Sigma = \{ u \in H^1_x(\R^d): xu \in L^2_x(\R^d) \}$, which appears frequently in the NLS theory, can now be interpreted naturally as the harmonic energy space; a scattering result for $\Sigma$ corresponds after the lens transformation to a local wellposedness result in the harmonic energy class for \eqref{lu-form} on the time interval $[-\pi/2,\pi/2]$.

Unfortunately, while the lens transform beautifully simplifies the pseudoconformal transformation, it makes the other symmetries of NLS, such as space or time translation symmetry $u(t,x) \mapsto u(t-t_0,x-x_0)$, 
the scaling symmetry \eqref{scaling}, and the Galilean invariance
\begin{equation}\label{gal}
 u(t,x) \mapsto e^{iv \cdot x} e^{-i|v|^2 t} u(t, x-vt)
 \end{equation}
somewhat more complicated (though still explicit, of course).  Only the compact symmetries, namely the phase rotation symmetry $u(t,x) \mapsto e^{i\theta} u(t,x)$, rotation symmetry $u(t,x) \mapsto u(t,U^{-1} x)$, and the time reversal 
symmetry $u(t,x) \mapsto \overline{u(-t,x)}$, remain unaffected by the lens transform.
Thus the lens transformation (which can be viewed as a quantization of the \emph{lens map} $(t,x) \mapsto (\tan t, \frac{x}{\cos t})$) ``straightens out'' the pseudoconformal transformation while distorting some of the other symmetries.

The lens transform also clarifies the relationship between the free Schr\"odinger equation 
$$(i \partial_t + \frac{1}{2} \Delta) u = 0$$ 
and  the linear Schr\"odinger equation with attractive harmonic potential
$$(i \partial_t + \frac{1}{2} \Delta - \frac{1}{2} |x|^2) v = 0,$$
as it maps solutions of the former to solutions of the latter without distorting initial data.  Indeed, since the former equation has the fundamental solution
$$ u(t,x) = e^{it \Delta/2} u_0(x) := \frac{1}{(2\pi i t)^{d/2}} \int_{\R^d} e^{i|x-y|^2/2t} u(0,y)\ dy$$
we obtain the \emph{Mehler formula} (see e.g. \cite{mehler}) for the fundamental solution of the latter\footnote{Strictly speaking, the lens transformation only permits one to verify the Mehler formula in the time interval $t \in (-\pi/2,\pi/2)$.  However, one then sees that solutions to this equation are periodic in time with period $2\pi$, either from this formula or via Hermite function expansion.  Of course, once a nonlinearity is introduced, there is no reason why the solution to the equation \eqref{lens} should continue to be time-periodic.}, namely
$$ v(t,x) = \frac{1}{(2\pi i \sin(t))^{d/2}} \int_{\R^d} e^{i|x-y|^2/2\tan(t)} v(0,y)\ dy.$$

The lens transform shows that the evolution of the NLS equation does not really ``stop'' at
time $t=+\infty$ or $t=-\infty$ (which correspond to the times $t=+\pi/2$ and $t=-\pi/2$ in the lens-transformed coordinates), but in fact continues on indefinitely beyond these points, providing of course that the equation
\eqref{lens} is globally wellposed; thus the apparent non-compactness of the time interval $\R$ is really an artefact of the coordinates rather than a true non-compactness of the flow (at least in the scattering-subcritical case $p > 1 + \frac{2}{d}$).  If the original solution $u: \R \times \R^d \to \C$ existed globally in time and enjoyed the asymptotic completeness relations
\begin{equation}\label{acr}
 \lim_{t \to \pm \infty} \| u(t) - e^{it\Delta/2} u_\pm \|_{L^2_x(\R^d)} = 0
\end{equation}
for some $u_\pm \in L^2_x(\R^d)$ (which we refer to as the \emph{asymptotic states} of $u$ at $t = \pm \infty$), 
then an easy computation using (recalling that $\L$ is linear and unitary on $L^2$ for any fixed time, and approximating $u_\pm$ in $L^2$ by Schwartz functions) shows that\footnote{Another way of viewing this is by observing that the Hermite propagator $e^{it(\Delta/2 - |x|^2/2)}$ (which quantizes a rotation of the phase plane $\{ (x,\xi): x, \xi \in \R^d \}$ by an angle $t$) is simply equal to the Fourier transform at $t = \pi/2$, which gives another explanation of the fact that this propagator is periodic of period $2\pi$.}
$$ \lim_{t \to \pm \pi/2} \| \L u(t) - \hat u_\pm \|_{L^2_x(\R^d)} = 0$$
where $\hat u_\pm$ is the Fourier transform of $u_\pm$:
$$ \hat u_\pm(x) := \int_{\R^d} e^{-ixy} u_\pm(y)\ dy.$$
Thus, asymptotic completeness in $L^2_x(\R^d)$ transforms under the lens transform to continuity in $L^2$ at the endpoint times $t=\pm \pi/2$; conversely, continuity of the lens-transformed solution at these times implies asymptotic completeness.  One can now cast the wave and scattering operators as the nonlinear propagators of \eqref{lens} between the times $t=0$, $t=+\pi/2$, and $t=-\pi/2$, composed with the Fourier or inverse Fourier transform at appropriate places; we omit the details.  This may conceptually clarify some of the scattering theory for NLS, such as that in \cite{caz}, \cite{ozawa}, \cite{blue}.  The long-time nonlinear propagator of \eqref{lens} can thus be factored into a product of the nonlinear propagator of \eqref{nls}, the scattering operator, and the Fourier transform.  The long-time evolution of \eqref{lens} is not well understood, even for small, spherically symmetric data in the defocusing case, precisely because it involves iterating the scattering operator (which can be viewed as a kind of Poincar\'e map for this flow), which is itself not well controlled.  See \cite{carles3} for further discussion of this iterated scattering operator.

\begin{remark} One can think of the lens time variable $\tan^{-1} t$ as an angular variable living on the unit circle $S^1$, or more precisely on the universal cover $\R$ of the circle.  Since the circle $S^1$ is the punctured plane $\R^2 \backslash \{0\}$ quotiented out by dilations, the lens time variable can thus be thought of as living in the universal cover of that punctured plane, quotiented out by dilations.  The ordinary time variable $t$ then lives on a vertical line $\{ (1,y): y \in \R \}$ of the punctured plane, lifted up to the universal cover and quotiented by dilations.  The pseudoconformal transformation \eqref{pcdef} $1/t$ for $t > 0$ uses a different coordinate chart, based on the horizontal line $\{ (x,1): x \in \R\}$ in the punctured plane, lifted up and quotiented as before.  (For $t < 0$ one uses the horizontal line $\{ (x,-1): x \in \R \}$.)  Thus the universal cover of the punctured plane, quotiented by dilations, can be viewed as a ``universal'' time axis, and all the other time variables viewed as co-ordinate charts.
\end{remark}

\begin{remark} There is an analogue of the lens transform for the repulsive harmonic potential, in which the trigonometric functions are replaced by their hyperbolic counterparts.  This transform is no longer a time compactification (it transforms $t$ to $\tanh^{-1} t$; instead, its \emph{inverse} is a time compactification, thus a bounded interval in the original time variable maps to an unbounded interval in the transformed variable.  See \cite{carles0}.
\end{remark}

\section{Application to pseudoconformal NLS}

To illustrate the simplifying conceptual power of the lens transform, we review the recent result in \cite{begout}, \cite{keraani} regarding equivalent formulations of the $L^2$ global existence problem for 
the pseudoconformal NLS 
\begin{equation}\label{nls-pc}
(i \partial_t + \frac{1}{2} \Delta) u = \mu |u|^{4/d} u.
\end{equation}
It is known (see \cite{gv}, \cite{cwI}, \cite{caz}) that for initial data $u_0 \in L^2_x(\R^d)$, there is a unique maximal Cauchy development $u: I \times \R^d \to \C$ in the space $C^0_{t,\loc} L^2_x(I \times \R^d) \cap L^{2(d+2)/d}_{t,\loc} L^{2(d+2)/d}_x(I \times \R^d)$ for some interval $I \subset \R$, and that this solution is global (so $I=\R$) if the mass $M(u) = M(u_0)$ is sufficiently small.  A similar statement is known for the lens-transformed equation
\begin{equation}\label{nls-pc-lens}
(i \partial_t + \frac{1}{2} \Delta - \frac{1}{2} |x|^2) v = \mu |v|^{4/d} v
\end{equation}
(see \cite{carles0}, \cite{carles}; one can also deduce this fact from the preceding one via the lens transform and time translation invariance).  We remark that the time-translation invariance $v(t,x) \mapsto v(t-t_0,x)$ of \eqref{nls-pc-lens} is manifestly obvious for the lens-transformed solution $v$, but corresponds to a rather unintuitive invariance for the original
solution $u$, namely
$$ u(t,x) \mapsto \frac{(1+s^2)^{d/4}}{(1+ts)^{d/2}} u( \frac{t-s}{1+ts}, \frac{x\sqrt{1+s^2}}{1+ts} ) e^{ i |x|^2 s / 2(1+ts) }$$
where $s := \tan t_0$; note that the pseudoconformal transformation \eqref{pcdef} is the limiting case $s=\pm \infty$.
The power of the lens transform lies in the fact that this non-trivial invariance of the original equation can be manipulated
effortlessly in the lens-transformed domain.

It is conjectured that the equation \eqref{nls-pc} has global solutions (with globally finite $L^{2(d+2)/d}_{t,x}$ norm) for any finite-mass initial data in the defocusing case $\mu=+1$, while in the focusing case $\mu=-1$, the same is conjectured for masses $M(u) < M(Q)$ less than that of the ground state $Q$, defined as the unique positive radial Schwartz solution to the ground state equation $\frac{1}{2} \Delta Q + Q^{1+\frac{4}{d}} = Q$.  Furthermore, the $L^{2(d+2)/d}$ norm is conjectured to be bounded by a quantity depending only on the mass.
As is well known, by starting with the soliton solution $u(t,x) = e^{it} Q(x)$ to \eqref{nls-pc} and applying the pseudoconformal transform one can obtain solutions of this equation with mass $M(Q)$ which blow up in finite time, which shows that this above conjecture is sharp.

As mentioned earlier, this conjecture is known to be true for sufficiently small mass, and has also recently been established for spherically symmetric data in a series of papers \cite{tvz}, \cite{tvz-2}, \cite{kvz}, \cite{ktv}.  At this time of writing, the conjecture remains open in general.

Ordinarily, global existence is a weaker statement than asymptotic completeness, or of global spacetime bounds.  However, a curious fact, recently established in \cite{keraani} in one and two dimensions, and in \cite{begout} in general dimensions (with the connection to \eqref{nls-pc-lens} essentially in \cite{carles0}), is
that these statements are in fact logically equivalent:

\begin{theorem}\label{main}\cite{keraani}, \cite{begout},\cite{carles0}  Let $p = 1 + \frac{4}{d}$, $\mu = \pm 1$, and $m > 0$ be fixed.  Then the following claims are equivalent.
\begin{itemize}
\item[(i)](Global wellposedness of \eqref{nls-pc}) For every $u_0 \in L^2_x(\R^d)$ with $M(u_0) < m$, the maximal Cauchy development $u$ to the equation \eqref{nls-pc} with data $u_0$ is global in time.
\item[(ii)](Asymptotic completeness of \eqref{nls-pc}) For every $u_0 \in L^2_x(\R^d)$ with $M(u_0) < m$, the maximal Cauchy development $u$ to the equation \eqref{nls-pc} with data $u_0$ is global in time, and furthermore there exist $u_+, u_- \in L^2_x(\R^d)$ for which \eqref{acr} holds.
\item[(iii)](Non-uniform spacetime bounds for \eqref{nls-pc}) For every $u_0 \in L^2_x(\R^d)$ with $M(u_0) < m$, the maximal Cauchy development $u$ to the equation \eqref{nls-pc} with data $u_0$ is global in time, and the $L^{2(d+2)/d}_{t,x}(\R \times \R^d)$ norm of $u$ is finite.
\item[(iv)](Uniform spacetime bounds for \eqref{nls-pc}) There exists a function $f: [0,m) \to [0,+\infty)$ such that for every $I \subset \R$ and $u \in C^0_t L^2_x(I \times \R^d) \cap L^{2(d+2)/d}_{t,x}(I \times \R^d)$ solving \eqref{nls-pc} with $M(u) < m$, we have the \emph{a priori} spacetime bound
$$ \| u \|_{L^{2(d+2)/d}_{t,x}(I \times \R^d)} \leq f(M(u)).$$
\item[(v)](Global wellposedness of \eqref{nls-pc-lens}) For every $u_0 \in L^2_x(\R^d)$ with $M(u_0) < m$, the maximal Cauchy development $u$ to the equation \eqref{nls-pc-lens} with data $u_0$ is global in time.
\item[(vi)](Uniform spacetime bounds for \eqref{nls-pc-lens}) There exists a function $f: [0,m) \to [0,+\infty)$ such that for every $I \subset \R$ and $v \in C^0_t L^2_x(I \times \R^d) \cap L^{2(d+2)/d}_{t,x}(I \times \R^d)$ solving \eqref{nls-pc-lens} with $M(v) < m$, we have the \emph{a priori} spacetime bound
$$ \| v \|_{L^{2(d+2)/d}_{t,x}(I \times \R^d)} \leq (1+|I|)^{d/2(d+2)} f(M(v)).$$
\end{itemize}
\end{theorem}

\begin{proof}

\textbf{Equivalence of $(i)$ and $(v)$}  From the lens transformation (which is mass-preserving, and also preserves the $C^0_t L^2_x$ and $L^{2(d+2)/d}_{t,x}$ norms) we see that global wellposedness of \eqref{nls-pc} is equivalent to local wellposedness of \eqref{nls-pc-lens} on the time interval $(-\pi/2,\pi/2)$.  But from time translation invariance (and mass conservation and uniqueness) we see that the local wellposedness of \eqref{nls-pc-lens} on this interval is equivalent to global wellposedness.

\textbf{Equivalence of $(v)$ and $(ii)$}  Clearly (ii) implies (i), which we know to be equivalent to (v).  On the other hand, (v) implies (ii) from the lens transform and the discussion in the introduction concerning \eqref{acr}.

\textbf{Equivalence of $(v)$ and $(iii)$}  Clearly (iii) implies (i), which we know to be equivalent to (v).  On the other hand, (v) implies (iii) since the lens transform preserves $L^{2(d+2)/d}_{t,x}$ norms.  Note how the local spacetime norm of $\L u$ is used to control the global spacetime norm of $u$.

\textbf{Equivalence of $(iv)$ and $(vi)$}  The implication of (vi) from (iv) follows from the fact that the lens transform preserves $L^{2(d+2)/d}_{t,x}$ norms.  Conversely, from (iv) and the lens transform we obtain (vi) for intervals $I$ contained in $(-\pi/2,\pi/2)$; the claim then follows from time translation invariance (and mass conservation and the triangle inequality).

\textbf{Implication of $(i)$ from $(iv)$} This follows immediately from the local wellposedness theory in \cite{cwI}, which among other things asserts that the maximal Cauchy development of an $L^2$ solution to \eqref{nls-pc} is not global only if the $L^{2(d+2)/d}_{t,x}$ norm is infinite.

\textbf{Implication of $(iv)$ from $(i)$} Note that we already know that (i) is equivalent to (ii), (iii), (v), so we may use these results freely.  This implication was established in \cite{begout}, following the induction-on-energy ideas of Keraani \cite{keraani} (see also \cite{mv}; related arguments also appear in \cite{borg:book}, \cite{gopher}, \cite{tvz-2}); in higher dimensions $d > 2$ the key harmonic analysis tool being the bilinear restriction estimate of the author \cite{tao}.  We can sketch a slightly simpler version of their argument (avoiding the full concentration-compactness machinery) as follows.  
Here we will avoid using the lens transform as it distorts the other symmetries $G$ of the NLS equation, which we will now need to exploit.

Let $\delta_0$ be the supremum of all $m$ for which (iv) holds; our task is to show that $\delta_0 \geq m$.  Suppose for contradiction that $\delta_0 < m$.  Then we could find
a sequence $u_n$ of global solutions to \eqref{nls-pc} with $M(u_n) < \delta_0$ such that
$\|u_n\|_{L^{2(d+2)/d}_{t,x}(\R \times \R^d)}$ was finite (by (iii)) but unbounded.  Let $u_{n,-}$ be the asymptotic scattering state at $t=-\infty$, thus $M(u_{n,-}) < \delta_0$ and $u_n(t)$ approaches $e^{it\Delta/2} u_{n,-}$ in $L^2_x(\R^d)$ as $t \to -\infty$.  From the small data theory we know that $\| e^{it\Delta/2} u_{n,-} \|_{L^{2(d+2)/d}_{t,x}(\R \times \R^d)} \geq \epsilon_d > 0$ for some absolute constant $\epsilon_d$ depending only on dimension, as otherwise $u_n$ would be bounded in $L^{2(d+2)/d}_{t,x}(\R \times \R^d)$ norm.
Applying the inverse Strichartz theorem (Theorem \ref{ist} in the appendix) in the contrapositive, we thus see that
$u_{n,-}$ is not weakly convergent with concentration to zero, or in other words there exist group elements $g_n \in G$ (where $G$ is defined in the appendix) such that $g_n u_{n,-}$ does not weakly converge to zero.  We observe that the group $G$ acts on $C^0_t L^2_x \cap L^{2(d+2)/d}_{t,x}$ solutions to \eqref{nls-pc} in a natural manner\footnote{In particular, the linear time translation operator \eqref{u0} acts on solutions to \eqref{nls-pc} by time translation (i.e. by the nonlinear propagator), while the modulation operator \eqref{mod} acts on solutions to \eqref{nls-pc} by the galilean invariance \eqref{gal}.} which is consistent with its action on the scattering data at $-\infty$, and which also preserves the mass and the $L^{2(d+2)/d}_{t,x}(\R \times \R^d)$ norm.  Thus without loss of generality we may take $g_n$ to be the identity for all $n$.  

Since
$u_{n,-}$ is bounded in $L^2_x(\R^d)$ and not weakly convergent to zero, we thus conclude by weak sequential compactness of the
$L^2_x$ ball that after passing to a subsequence (which we continue to call $u_{n,-}$), there exists a non-zero $u_- \in L^2_x(\R^d)$ with $0 < M(u_-) \leq \delta_0 < m$ such that $u_{n,-}$ converges weakly to $u_-$.  By the hypothesis (v) and the lens transform, we can find a global solution $u \in C^0_t L^2_x \cap L^{2(d+2)/d}_{t,x}$ with finite $L^{2(d+2)/d}_{t,x}$ norm
which has $u_-$ as its asymptotic state at $-\infty$.  Also, if we split $u_{n,-} = u_- + v_{n,-}$, then from Pythagoras's theorem and the weak convergence of $v_{n,-}$ to zero we have 
$$\limsup_{n \to \infty} M(v_{n,-}) = \limsup_{n \to \infty} M(u_{n,-}) - M(u_-) \leq \delta_0 - M(u_-) < \delta_0.$$
Thus by construction of $\delta_0$, we thus see that (for $n$ sufficiently large) there are global solutions $v_n \in C^0_t L^2_x \cap L^{2(d+2)/d}_{t,x}$ to \eqref{nls-pc} which have $v_{n,-}$ has their asymptotic state, and whose $L^{2(d+2)/d}_{t,x}$ norms are uniformly bounded in $n$.  Since $v_{n,-}$ converges weakly to zero, it is not difficult (using the Strichartz wellposedness theory and the uniform spacetime bounds) to show that $v_n$ also converges weakly to zero\footnote{This is easiest to establish by duality, viewing $v_n$ as the linear evolution of the variable coefficient (but self-adjoint) Schr\"odinger operator $i \partial_t + \frac{1}{2} \Delta - \mu |v_n|^{4/d}$ with asymptotic state $v_{n,-}$ at $t=-\infty$.  Testing $v_n$ against a spacetime test function $\varphi$ then reduces one to studying the solution $w$ of the variable coefficient inhomogeneous equation $(i \partial_t + \frac{1}{2} \Delta - \mu |v_n|^{4/d}) w = \varphi$ with asymptotic state $0$ at $t=+\infty$, but this can be controlled by Strichartz estimates.}.  

Now split
$$ u_n = u + v_n + e.$$
One checks that $e(t) \to 0$ in $L^2_x(\R^d)$ norm as $t \to -\infty$, and that $e$ solves the equation
$$ (i \partial_t + \Delta) e = [F(u+v_n+e)-F(u+v_n)] + [F(u+v_n) - F(u)-F(v_n)]$$
where $F(z) := \mu |z|^{4/d} z$.  Since $u,v_n$ are uniformly bounded in $L^{2(d+2)/d}_{t,x}(\R \times \R^d)$ and $v_n$ converges weakly to zero, one easily verifies that $F(u+v_n) - F(u)-F(v_n)$ converges strongly in $L^{2(d+2)/(d+4)}_{t,x}(\R \times \R^d)$ to zero (cf. \cite[Lemma 5.5]{begout}).  Standard application of Strichartz wellposedness theory (estimating
$F(u+v_n+e)-F(u+v_n)$ pointwise by $O(|e| ( |e| + |u| + |v_n| )^{4/d})$; see e.g. \cite{tvz-2})
and the uniform bounds on $u,v_n$ then shows that $e$ converges strongly in $L^{2(d+2)/d}_{t,x}(\R \times \R^d)$ to zero.  But this then implies that $u_n$ is uniformly bounded in $L^{2(d+2)/d}_{t,x}(\R \times \R^d)$, a contradiction.
\end{proof}

\begin{remark} In one dimension $d=1$, the equivalent statements in Theorem \ref{main} are also linked to the analogous statements for the $L^2$-critical generalised Korteweg-de Vries equation, see \cite{tao-gkdv}.  These statements are also stable under addition of further power nonlinearities which are larger than the $L^2$-critical power $1+\frac{4}{d}$, but less than or equal to the $H^1$-critical power $1 + \frac{4}{d-2}$; see \cite{zhang}.
\end{remark}

\appendix

\section{The inverse Strichartz theorem}

The standard Strichartz estimate (see e.g. \cite{strich}) asserts that
\begin{equation}\label{itdelta}
\| e^{it\Delta/2} u_0 \|_{L^{2(d+2)/d}_{t,x}(\R \times \R^d)} \leq C_d \| u_0 \|_{L^2_x(\R^d)}
\end{equation}
for all $u_0 \in L^2_x(\R^d)$ and some constant $0 < C_d < \infty$.  It is of interest to invert this estimate by
deducing necessary and sufficient conditions for which this estimate is sharp.  Such an inverse result is implicitly
in \cite{begout} (following \cite{keraani} and \cite{mv}), but we state it explicitly here.

We first observe that both sides of \eqref{itdelta} are invariant under the spatial translation symmetry
$$ u_0(x) \mapsto u_0(x - x_0)$$
for any $x_0 \in \R$, the phase rotation symmetry
$$ u_0(x) \mapsto e^{i\theta} u_0(x)$$
for any $\theta \in \R$, the scaling symmetry
$$ u_0(x) \mapsto \frac{1}{\lambda^{d/2}} u_0(\frac{x}{\lambda})$$
for any $\lambda > 0$, the time translation symmetry
\begin{equation}\label{u0}
u_0(x) \mapsto e^{-it_0\Delta/2} u_0
\end{equation}
for any $t_0 \in \R$, and the modulation symmetry
\begin{equation}\label{mod}
u_0(x) \mapsto e^{i v \cdot x} u_0
\end{equation}
for any $v \in \R^d$.  (The latter is not immediately obvious, but follows from the Galilean invariance \eqref{gal}
of the free linear Schr\"odinger equation $iu_t + \frac{1}{2} \Delta u = 0$.)  Let us use $G$ to refer to the group of unitary transformations on $L^2_x(\R^d)$ generated by all these symmetries; this is a $2d+3$-dimensional non-compact Lie group whose elements can be explicitly described, though we will not do so here\footnote{The estimate \eqref{itdelta} is also invariant under the pseudoconformal transformation, rotation symmetry, the Fourier transform, and the quadratic modulation symmetry $u_0(x) \mapsto e^{i\alpha |x|^2} u_0(x)$, which increases the dimension of the symmetry group to $\frac{d^2}{2} + \frac{3d}{2} + 4$, but we will not use these symmetries here as they are not needed for the inverse Strichartz theorem.}.  Thus for all $g \in G$ and $u_0 \in L^2_x(\R^d)$ we have
$$ \| e^{it\Delta/2} g u_0 \|_{L^{2(d+2)/d}_{t,x}(\R \times \R^d)}
= \| e^{it\Delta/2} u_0 \|_{L^{2(d+2)/d}_{t,x}(\R \times \R^d)}$$
and
$$ \| g u_0 \|_{L^2_x(\R^d)} = \| u_0 \|_{L^2_x(\R^d)}.$$

We remark that this group $G$ is a group of dislocations in the sense of Schindler and Tintarev \cite{tintarev}.  In other words, if $g_n \in G$ is any sequence of group elements in $G$ which has no strongly convergent subsequence (in the strong operator topology on $L^2_x(\R^d)$), then $g_n$ necessarily converges to zero in the weak operator topology.

Let us say that a bounded sequence $u_n \in L^2_x(\R^d)$ \emph{converges weakly with concentration} to zero if the sequence $g_n u_n$ is weakly convergent to zero for every choice of group element $g_n \in G$, thus $\lim_{n \to \infty} \langle g_n u_n, v \rangle_{L^2_x(\R^d)} = 0$ for all $g_n \in G$ and $v \in L^2_x(\R^d)$.  This convergence is a little stronger than weak convergence, but certainly weaker than strong convergence.  Nevertheless, for the purposes of Strichartz estimates, it is ``as 
good as'' strong convergence in the following sense:

\begin{theorem}[Inverse Strichartz theorem]\label{ist}  Suppose that $u_n \in L^2_x(\R^d)$ is a bounded sequence in $L^2_x(\R^d)$ which converges weakly with concentration to zero.  Then
\begin{equation}\label{sw}
\lim_{n \to \infty} \| e^{it\Delta/2} u_n \|_{L^{2(d+2)/d}_{t,x}(\R \times \R^d)} = 0.
\end{equation}
\end{theorem}

\begin{proof}  Fix the sequence $\vec u = (u_n)_{n =1}^\infty$.  We shall use the following asymptotic notation:
\begin{itemize}
\item We use $O(X)$ to denote any quantity bounded in magnitude by $C(\vec u,d) X$, where $0 < C(\vec u,d) < \infty$ depends only on $\vec u$ and $d$.
\item We use $o(X)$ to denote any quantity bounded in magnitude by $c(n,\vec u,d) X$, where $\lim_{n \to \infty} c(n,\vec u,d) = 0$ for each fixed $\vec u$ and $d$.  
\end{itemize}
Thus for instance our hypotheses on $\vec u$ imply that 
$$ \|u_n\|_{L^2_x(\R^d)} = O(1) \hbox{ and } \langle u_n, g_n v \rangle = o_v(1)$$
for each fixed $v \in L^2_x(\R^d)$ and all $g_n \in G$, where we subscript $o(1)$ by $v$ to indicate that the implied constant $c(n,\vec u,d)$ can depend on $v$.  Our objective is to show that $\| e^{it\Delta/2} u_n \|_{L^{\frac{2(d+2)}{d}}_{t,x}(\R \times \R^d)} = o(1)$.

We consider the mesh of dyadic cubes $Q$ in $\R^d$ (i.e. half-open cubes $Q$ with axes parallel to the coordinate axes, whose length $\ell(Q)$ is a power of two, and whose corners have coordinates which are integer multiples of the length).  These cubes should be thought of as lying in Fourier space rather than in physical space.
We say that two cubes $Q,Q'$ are \emph{close} if they have the same length and are not adjacent (i.e. their closures do not intersect), but their parents are adjacent.  For any cube $Q$, let $u_{n,Q}$ be the Fourier restriction of $u_n$ to $Q$, thus $\hat u_{n,Q} = 1_Q \hat u_n$.  We write $u_{n,Q}(t,x)$ for the free solution $e^{it\Delta/2} u_{n,Q}(x)$.

Observe that given any two distinct frequencies $\xi,\xi' \in \R^d$ there is a unique pair of close dyadic cubes $Q,Q'$ which contain $\xi,\xi'$ respectively.  This gives rise (as in \cite{tvv}) to the Whitney decomposition
$$ u_n^2 = \sum_{Q,Q' \hbox{ close}} u_{n,Q} u_{n,Q'}$$
and hence
$$ \| u_n \|_{L^{\frac{2(d+2)}{d}}_{t,x}(\R \times \R^d)}^2
= \| \sum_{Q,Q' \hbox{ close}} u_{n,Q} u_{n,Q'} \|_{L^{\frac{d+2}{d}}_{t,x}(\R \times \R^d)}.$$
Thus it will suffice to show that
\begin{equation}\label{targ}
\| \sum_{Q,Q' \hbox{ close}} u_{n,Q} u_{n,Q'} \|_{L^{\frac{d+2}{d}}_{t,x}(\R \times \R^d)} = o(1).
\end{equation}

Since $u_n$ converges weakly to zero, one can easily show that $u_{n,Q}(t,x)$ converges pointwise to zero for each fixed $t,x,Q$.  However, this convergence is not uniform in $t,x,Q$.  By using the stronger hypothesis that $u_n$ converges weakly \emph{with concentration} to zero, however, one obtains the more uniform estimate
\begin{equation}\label{unq}
 u_{n,Q}(t,x) = o( \ell(Q)^{d/2} )
\end{equation}
for all $t,x,Q$.  Indeed one can use the group $G$ to move $t,x,Q$ to a compact set, in which case the uniform estimate follows from the pointwise estimate (observe that the frequency localisation to $Q$ and the uniform $L^2$ bounds on $u_n$ ensure that the $u_{n,Q}$ are equicontinuous in $n$).  

Let $Q,Q'$ be two close cubes.  The bilinear restriction theorem from \cite{tao}, combined with a standard parabolic rescaling argument (see e.g. \cite{tvv} or \cite[Corollary 2.3]{begout}) then ensures that
$$ \| u_{n,Q} u_{n,Q'} \|_{L^q_{t,x}(\R \times \R^d)} = O(\ell(Q)^{d - \frac{d+2}{q}})$$
for any fixed $\frac{d+2}{d} > q > \frac{d+3}{d+1}$.  On the other hand, from \eqref{unq} we have
$$ \| u_{n,Q} u_{n,Q'} \|_{L^\infty_{t,x}(\R \times \R^d)} = o(\ell(Q)^d)$$
and hence by interpolation we have
\begin{equation}\label{unqq}
 \| u_{n,Q} u_{n,Q'} \|_{L^{\frac{d+2}{d}}_{t,x}(\R \times \R^d)} = o(1).
\end{equation}
The point is that the right-hand side is uniform over all choices of close cubes $Q,Q'$.  

Now we sum over all close cubes $Q,Q'$.  Observe (as in \cite{tvv}, \cite{begout}, or in the classical work of C\'ordoba \cite{cordoba}) that the spacetime Fourier transforms of $u_{n,Q} u_{n,Q'}$ are supported in essentially disjoint cubes.  Thus we may apply the almost orthogonality
estimate in \cite[Lemma 6.1]{tvv} to conclude
$$ \| \sum_{Q,Q' \hbox{ close}} u_{n,Q} u_{n,Q'} \|_{L^{\frac{d+2}{d}}_{t,x}(\R \times \R^d)} 
= O( (\sum_{Q,Q' \hbox{ close}} \| u_{n,Q} u_{n,Q'} \|_{L^{\frac{d+2}{d}}_{t,x}(\R \times \R^d)}^p)^{1/p} )$$
for some $1 < p < \infty$ depending only on $d$.  But from the bilinear restriction estimate from \cite[Theorem 2.3]{tvv} 
(or by interpolating the estimate in \cite{tao} with trivial estimates) and a parabolic rescaling argument we have
$$ \| u_{n,Q} u_{n,Q'} \|_{L^{\frac{d+2}{d}}_{t,x}(\R \times \R^d)} = 
O( |Q|^{1 - \frac{2}{q}} \|\hat u_{n}(0)\|_{L^q_\xi(Q)} \|\hat u_{n}(0)\|_{L^q_\xi(Q')} )
$$
for some $1 < q < 2$ whose exact value is not important here.

Suppose for the moment that we could show
\begin{equation}\label{show}
(\sum_{Q,Q' \hbox{ close}} [|Q|^{1 - \frac{2}{q}} \|\hat u_{n}(0)\|_{L^q_\xi(Q)} \|\hat u_{n}(0)\|_{L^q_\xi(Q')}]^r)^{1/r} = O_r(\|\hat u_n(0)\|_{L^2_\xi(\R^d)}^2)
\end{equation}
for all $r > 1$.  Then by the uniform boundedness of $u_n$ in $L^2_x$ we would have
$$ (\sum_{Q,Q' \hbox{ close}} \| u_{n,Q} u_{n,Q'} \|_{L^{\frac{d+2}{d}}_{t,x}(\R \times \R^d)}^r)^{1/r} = O_r(1)$$
and hence by \eqref{unqq} and by choosing $1 < r < p$
$$ (\sum_{Q,Q' \hbox{ close}} \| u_{n,Q} u_{n,Q'} \|_{L^{\frac{d+2}{d}}_{t,x}(\R \times \R^d)}^p)^{1/p} = o(1)$$
from which \eqref{targ} follows.  Thus it suffices to show \eqref{show}.  Using the elementary inequality $ab \leq \frac{1}{2} a^2 + \frac{1}{2} b^2$ and noting that each cube $Q$ has only $O(1)$ cubes $Q'$ that are close to it, we reduce to showing that
$$ (\sum_Q [|Q|^{1 - \frac{2}{q}} \|\hat u_{n}(0)\|_{L^q_\xi(Q)}^2]^r)^{1/r} = O(\|\hat u_n(0)\|_{L^2_\xi(\R^d)}^2).$$
This follows from \cite[Theorem 1.3]{begout}; for the convenience of the reader we give a short proof here.  Setting
$f := |\hat u_n(0)|^q$, $p := 2/q$, and $s := rp$ it suffices to show that
$$ (\sum_Q |Q|^{\frac{s}{p} - s} \|f\|_{L^1_\xi(Q)}^s)^{1/s} = O_{s,p}( \|f\|_{L^p_\xi(\R^d)} )$$
for all $s > p > 1$ and $f \in L^p_\xi(\R^d)$. By the real interpolation method it suffices to prove the restricted estimate
$$ (\sum_Q |Q|^{\frac{s}{p} - s} |\Omega \cap Q|^s)^{1/s} = O_{s,p}( |\Omega|^{1/p} )$$
for all $s > p > 1$ and all sets $\Omega$ of finite measure.  But since $|\Omega \cap Q|^s \leq |\Omega \cap Q| \min(|\Omega|,|Q|)^{s-1}$ we reduce to showing that
$$ \sum_k \sum_{Q: \ell(Q) = k} |Q|^{\frac{s}{p}-s} |\Omega \cap Q| \min(|\Omega|,|Q|)^{s-1} = O_{s,p}(|\Omega|^{s/p}).$$
But the left-hand side sums to
$$ \sum_k (2^{dk})^{\frac{s}{p}-s} |\Omega| \min(|\Omega|,2^{dk})^{s-1}$$
which can be computed to be $O_{s,p}(|\Omega|^{s/p})$ as claimed.
\end{proof}

\begin{remark}  Theorem \ref{ist} also follows immediately from the concentration compactness theorem in \cite[Theorem 5.4]{begout}; conversely, that theorem follows quickly from Theorem \ref{ist}, \cite[Lemma 5.5]{begout}, and the abstract concentration compactness theorem in \cite{tintarev}.  However the proof above, while using many of the same ingredients as that in \cite{begout}, uses slightly less machinery and thus can be regarded as a more primitive (but non-quantitative) proof.
\end{remark}

One can state Theorem \ref{ist} in the contrapositive, in a manner which more clearly explains the terminology ``inverse Strichartz theorem'':

\begin{corollary}[Inverse Strichartz theorem, again]\label{ist-cor}  Let $d \geq 1$ and $m, \eps > 0$.  Then there exists a finite set ${\mathcal C} \subset L^2_x(\R^d)$ of functions of norm $1$ and an $\eta > 0$ with the following property: whenever $u \in L^2_x(\R^d)$ obeys the bounds
\begin{equation}\label{weakconc}
\| u\|_{L^2_x(\R^d)} \leq m; \quad \| e^{it\Delta/2} u \|_{L^{2(d+2)/d}_{t,x}(\R \times \R^d)} \geq \eps
\end{equation}
then there exists $g \in G$ and $v \in {\mathcal C}$ such that $|\langle u, g v \rangle_{L^2_x(\R^d)}| \geq \eta$.
\end{corollary}

\begin{proof} Assume for contradiction that the corollary failed, then there exist $d,m,\eps$ and a sequence $u_n$, each of which obeys \eqref{weakconc}, which is weakly convergent with concentration to zero.  But this contradicts Theorem \ref{ist}.
\end{proof}

\begin{remark} By going more carefully through the arguments in \cite{begout}, one can obtain more quantitative estimates here; indeed, for any fixed $d$, the quantities $\eta$ and $\# {\mathcal C}$ will be some polynomial combination of $m$ and $\eps$.  Furthermore, for any fixed $k \geq 1$, we can make the elements of ${\mathcal C}$ bounded in the weighted Sobolev space $H^{k,k}_x(\R^d)$ with a norm which is polynomial in $m$ and $\eps$.  We omit the details.
\end{remark}


\begin{thebibliography}{10}

\bibitem{begout}
P. Begout, A. Vargas, \emph{Mass concentration phenomena for the $L^2$-critical nonlinear Schr\"odinger equation},  Trans. Amer. Math. Soc.  \textbf{359}  (2007),  no. 11, 5257--5282.

\bibitem{blue}
P. Blue, J. Colliander, \emph{Global well-posedness in Sobolev space implies global existence for weighted $L^2$ initial data for $L^2$-critical NLS}, preprint.

\bibitem{borg:book}
J. Bourgain, \emph{New global well-posedness results for non-linear Schr\"odinger equations}, AMS Publications, 1999.

\bibitem{carles0}
R. Carles, \emph{Critical nonlinear Schr\"odinger equations with and without harmonic potential}, Math. Models Methods Appl. Sci. \textbf{12} (2002), no. 10, 1513--1523.

\bibitem{carles}
R. Carles, \emph{Nonlinear Schr\"odinger equation with harmonic potential and applications}, SIAM J. Math. Anal. \textbf{35} (2003), 823--843.

\bibitem{carles3}
R. Carles, \emph{Semi-classical Schr\"odinger equations with harmonic potential and nonlinear perturbation}, Ann. Inst. H. Poincar\'e Anal. Nonlin\'eaire \textbf{20} (2003), no. 3, 501--542.

\bibitem{carles2}
R. Carles, \emph{Linear vs. nonlinear effects for nonlinear Schr\"odinger equation with potential}, Commun. Contemp. Math. \textbf{7} (2005), no. 4, 483--508.

\bibitem{caz}
T. Cazenave, \emph{Semilinear Schr\"odinger equations}, Courant
Lecture Notes in Mathematics, \textbf{10}. New York University,
Courant Institute of Mathematical Sciences, AMS, 2003.

\bibitem{cwI}
T. Cazenave, F.B. Weissler, \emph{Critical nonlinear Schr\"odinger
Equation}, Non. Anal. TMA, \textbf{14} (1990), 807--836.

\bibitem{cw-rap}
T. Cazenave, F.B. Weissler, \emph{Rapidly decreasing solutions of the nonlinear Schr\"odinger equation}, Comm. Math. Phys. \textbf{147} (1992), no. 1, 75--100.

\bibitem{christ}
D. Christodoulou, S. Klainerman, \emph{Asymptotic properties of linear field equations in Minkowski space}, Comm. Pure Appl. Math. \textbf{43} (1990), 137--199.

\bibitem{gopher}
J. Colliander, M. Keel, G. Staffilani, H. Takaoka, T. Tao, \emph{Global well-posedness and scattering in the energy space for the critical nonlinear Schrodinger equation in $\R^3$}, to appear, Annals Math.

\bibitem{cordoba}
A. C\'ordoba, \emph{The Kakeya maximal function and the spherical summation multipliers}, Amer. J. Math. \textbf{99} (1977), 1--22.

\bibitem{mehler}
R.P. Feynman, A.R. Hibbs, Quantum mechanics and path integrals (International Series in Pure and Applied Physics), Maidenhead, Berksh.: McGraw-Hill Publishing Company Ltd. 365 p., 1965.

\bibitem{gv}
J. Ginibre, G. Velo, \emph{On a class of nonlinear Schr\"odinger problems I. The Cauchy problem, general case}, J. Funct. Anal. \textbf{32} (1979), 1--32.

\bibitem{tao:keel}
M. Keel, T. Tao, \emph{Endpoint Strichartz Estimates}, Amer. Math. J. 
120 (1998), 955--980.

\bibitem{keraani-old}
S. Keraani, \emph{On the Defect of Compactness for the Strichartz Estimates of the Schr\"odinger Equations}, J. Differential Equations \textbf{175} (2001), no. 2, 353--392.

\bibitem{keraani}
S. Keraani, \emph{On the blow-up phenomenon of the critical nonlinear Schr\"odinger equation},  J. Funct. Anal.  \textbf{235}  (2006),  no. 1, 171--192.

\bibitem{ktv}
R. Killip, T. Tao, M. Visan, \emph{The cubic nonlinear Schrodinger equation in two dimensions with radial data}, preprint.

\bibitem{kvz}
R. Killip, M. Visan, X. Zhang, \emph{The mass-critical nonlinear Schrodinger equation with radial data in dimensions three and higher}, Analysis and PDE \textbf{1} (2008), 229--266.

\bibitem{mv}
F. Merle, L. Vega, \emph{Compactness at blow-up time for $L^2$ solutions of the critical nonlinear Schr\"odinger equation in 2D}, Internat. Math. Res. Not. \textbf{8} (1998), 399--425.

\bibitem{lens}
U. Niederer, \emph{The maximal kinematical invariance groups of Schr\"odinger equations with arbitrary potentials}, Helv. Phys. Acta \textbf{47} (1974), 167--172.

\bibitem{ozawa}
T. Ozawa, \emph{Long range scattering for nonlinear Schrodinger equations in one space dimension}, Commun. Math. Phys. 139 (1991), 479--493.

\bibitem{penrose}
R. Penrose, \emph{Conformal treatment of infinity}, 1964 Relativit\'e, Groupes et Topologie, pp. 565--584, Gordon and Breach, New York.

\bibitem{lens2}
A.V. Rybin, G.G. Varzugin, M. Lindberg, J. Timonen, R.K. Bullough, \emph{Similarity solutions and collapse in the attractive Gross-Pitaevskii equation}, Phys. Rev. E. \textbf{62} (2000), no. 5, part A, 6224--6228.

\bibitem{tintarev}
I. Schinder, K. Tintarev, \emph{An abstract version of the concentration compactness principle}, Revista Mathem\'atica Complutense \textbf{15} (2002), 417--436.

\bibitem{strich}
R.~S. Strichartz, \emph{Restriction of Fourier Transform to Quadratic Surfaces
and Decay of Solutions of Wave Equations}, Duke Math. J., \textbf{44} (1977), 70
5--774.

\bibitem{tao}
T. Tao, \emph{A sharp bilinear restriction estimate for paraboloids}, Geom. Func. Anal. \textbf{13} (2003), 1359--1384.

\bibitem{tao-gkdv}
T. Tao, \emph{Two remarks on the generalised Korteweg-de Vries equation}, Discrete Cont. Dynam. Systems \textbf{18} (2007), 1--14.

\bibitem{tvv}
T. Tao, A. Vargas, L. Vega, \emph{A bilinear approach to the restriction and Kakeya conjectures}, J. Amer. Math. Soc. \textbf{11} (1998), 967--1000.

\bibitem{zhang}
T. Tao, M. Visan, X. Zhang, \emph{The nonlinear Schr\"odinger equation with combined power-type nonlinearities}, Comm. PDE \textbf{32} (2007), 1281--1343.

\bibitem{tvz}
T. Tao, M. Visan, X. Zhang, \emph{Global well-posedness and scattering for the mass-critical nonlinear Schrodinger equation for radial data in high dimensions}, Duke Math. J. \textbf{140} (2007), 165--202.

\bibitem{tvz-2}
T. Tao, M. Visan, X. Zhang, \emph{Minimal-mass blowup solutions of the mass-critical NLS}, Forum Mathematicum  \textbf{20} (2008), 881--919.

\end{thebibliography}
\end{document}